\documentclass[12pt]{article}

\usepackage{amsmath}
\usepackage{amssymb}
\usepackage{amsfonts}
\usepackage{graphics}

\newtheorem{thm}{Theorem}[section]
\newtheorem{lemma}[thm]{Lemma}
\newtheorem{cor}[thm]{Corollary}

\newtheorem{prop}[thm]{Proposition}

\newenvironment{remark}{\par\medskip\noindent{\bf Remark.\ }}{\par\smallskip}

\newcommand{\be}{\begin{equation}}
\newcommand{\ee}{\end{equation}}
\newcommand{\openbox}{\leavevmode
  \hbox to8pt{\hfil\vrule\vbox to6pt{\hrule width6pt\vfil\hrule}\vrule}}

\newcommand{\qed}{\hbox to5pt{ } \hfill \openbox\bigskip\medskip}

\newcommand{\func}{{\mbox{{\sf func}}}}

\newcommand{\Fq}{\mathbb F _q}
\newcommand{\Zp}{\mathbb Z _p}

\newcommand{\cV}{\mbox{$\cal V$}}

\newcommand{\cF}{\mbox{$\cal F$}}
\newcommand{\cG}{\mbox{$\cal G$}}
\newcommand{\cH}{\mbox{$\cal H$}}

\newcommand{\AC}{\mbox{$AC$}}
\newcommand{\ac}{\mbox{$ac$}}
\newcommand{\ve}[1]{\mathbf{#1}}

\newcommand{\N}{\mathbb N}
\newcommand{\Z}{\mathbb Z}
\newcommand{\Q}{\mathbb Q}

\newcommand{\F}{\mathbb F}

\newcommand{\Sm}{\mbox{\rm Sm}}

\newcommand{\Sym}{\mathop\textup{Sym}}
\newcommand{\Aff}{\mathop\textup{Aff}}
\title{Almost covers of finite sets of points}
\author{G\'abor Heged\"{u}s
\\{\normalsize  \'Obuda University}
\\{\normalsize B\'ecsi \'ut 96/B, Budapest, Hungary, H-1032}
\\{\normalsize hegedus.gabor@uni-obuda.hu}
}

\begin{document}

\maketitle
\begin{abstract}
Let
$\mbox{$\cal V$} \subseteq {\mathbb F}^n$ be a finite set of points in an affine space. A finite set of affine hyperplanes $\{H_1, \ldots ,H_m\}$ is said to be an {\em almost cover} of $\mbox{$\cal V$}$ and $\mathbf{v}$, if their union $\cup_{j=1}^m H_j$ contains
$\mbox{$\cal V$}\setminus \{\mathbf{v}\}$ but does not contain $\mathbf{v}$. 

We give here a lower bound for the size of a minimal almost cover of $\mbox{$\cal V$}$ and $\mathbf{v}$ in terms of the size of $\mbox{$\cal V$}$ and the dimension $n$. 

We prove a  generalization of  Sziklai  and  Weiner's Theorem. 
Our simple proof is based on Gr\"obner basis theory.
\end{abstract}
\medskip

\section{Introduction}

Throughout this paper $n$ denotes a positive integer, and $[n]$ stands for 
the set $\{1,2,\dots, n\}$. We denote 
by $2^{[n]}$ the set of all subsets of $[n]$. Subsets of $2^{[n]}$ are called \emph{set families}.
Let $\binom{[n]}m$ stand for the family of all subsets of $[n]$
which have cardinality $m$. Let $\binom{[n]}{\le m}$ denote the family  of all subsets that 
have size at most $m$. 

Let  $\F$ be a field and let 
$\F[x_1, \ldots, x_n]=\F[\ve x]$ denote  the
ring of polynomials in the variables $x_1, \ldots, x_n$ over $\F$.
For a subset $F \subseteq [n]$ we write
$x_F = \prod_{j \in F} x_j$.
In particular, $x_{\emptyset}= 1$. We denote by $\ve v_F\in \{0,1\}^n$ the
characteristic vector of a set
$F \subseteq [n]$.
 
Define
$V(\cF)$ as the subset  $\{\ve v_F : F \in \cF\} \subseteq \{0,1\}^n \subseteq \F^n$ for any family of subsets $\cF \subseteq 2^{[n]}$.

It is natural to consider the vanishing ideal $I(V(\cF))$:
$$ 
I(V(\cF)):=\{f\in \F[\ve x]:~f(\ve v)=0 \mbox{ whenever } \ve v\in V(\cF)\}. 
$$
It is easy to verify that we can identify the algebra $\F[\ve x]/I(V(\cF))$ and the algebra of
$\F$ valued functions on $V(\cF)$. Consequently 
$$
\dim _\F
\F[\ve x]/I(V(\cF))=|\cF|.
$$

Let $0\leq k<n$ be integers. Consider the set
of vectors
$$
\cV(n,k):=\{\ve v_K: K\in{[n]\choose {\leq	k}}\}\subseteq \{0,1\}^n \subseteq \F^n.
$$
Let $T\subseteq [n]$ be a fixed subset such that $|T|>k$. Define
$$
\cV(n,k,T):=\cV(n,k)\cup \{\ve v_T\}.
$$

Let $q>1$ be a positive integer.
Suppose that $|\F|\geq q$. Then we can consider $[q]$ as a subset of
our ground field $\F$ via an injective map $i: [q]\rightarrow \F$.  
We can consider the set
of non-decreasing sequences
$$
I(n,q)=\{(f_1,f_2, \dots, f_n) \in [q]^n:~ f_1\leq f_2\leq\cdots \leq f_n \}.
$$
Clearly we have
$$
|I(n,q)|={n+q-1 \choose q-1}.
$$

Denote
by $J(n,q)$ the image of $I(n,q)$ by the map induced by $i$:
$$ 
J(n,q)=\{ (i(v_1), \ldots ,i(v_n)): ~~ (v_1,\ldots v_n)\in
I(n,q) \}.
$$

Let
$\cV \subseteq \F^n$ be a finite set of points in an affine space. Then a finite set of affine hyperplanes $\{H_1, \ldots ,H_m\}$ is said to be an {\em almost cover} of $\cV$ and $\ve v$, if their union $\cup_{j=1}^m H_j$ contains
$\cV\setminus \{\ve v\}$ but does not contain $\ve v$.

Jamison proved  the following  classical result in \cite{J}: an almost cover of the
$n$-dimensional affine space over the $q$-element finite field requires at
least $(q-1)n$ hyperplanes. 
A. Blokhuis, A.E. Brouwer, and T. Sz\H onyi proved further results in finite geometries  in \cite{BBS}.

J. Aaronson, C. Groenland, A. Grzesik, T. Johnston and B. Kielak
defined in \cite{AGGJK} the {\em exact cover} of any subset $\cV \subseteq \{0,1\}^n$ as follows: the exact cover of $\cV$ is a set of hyperplanes whose union intersect $\{0,1\}^n$ exactly in $\cV$, i.e. the affine points of  $\{0,1\}^n\setminus \cV$ is not covered. Following \cite{AGGJK} we denote the exact cover number of $\cV$ by $\mbox{ec}(\cV)$:  here  $\mbox{ec}(\cV)$ is the minimum number of an exact cover of $\cV$. 
Consequently $\mbox{ec}(\cV(n,k))=k+1$.
Brouwer and Schrijver proved an equivalent version in \cite{BS2}: to pierce every affine hyperplane in
$\F_q^n$ one needs at least $(q-1)n+1$ points. 

We proved the following result about the almost cover of the set $J(n,q)$ in \cite{HR} Theorem 1.2: 
Let $0<k\leq n$ be an integer and let $\ve s_1,\ldots ,\ve s_k\in J(n,q)$ be
increasing vectors.
Let $\{H_j:~ 1\leq j\leq m\}$ be a set of affine hyperplanes  such that
$$
J(n,q)\setminus \{\ve s_1, \ldots ,\ve s_k\}\subseteq \cup_{j=1}^m H_j.
$$
Then $m\geq q-1$.

Our work is motivated by the
Alon-F\"uredi theorem for the hypercube \cite{AF}: {\em Every almost cover of the vertex set
of an $n$-dimensional cube requires at least $n$ affine hyperplanes.}

R. Karasev proved the following version of  the Cayley-Bacharach Theorem in \cite{K}.

\begin{thm}\label{Kar}
Let  $\F$ be a field and $g_1, \ldots ,g_m\in \F[x_1,x_2,\ldots, x_n]$ have degrees $k_1,\ldots,k_n$. 
Assume the system $g_1=\cdots=g_n=0$ has exactly $K=k_1\cdots k_n$ isolated solutions in $\mathbb{P}^{n}$. Then for every polynomial $f$ with  $\deg f\le (k_1+\cdots+k_n)-n-1$, there exists a nontrivial linear relation
$$
\sum_{x\in X}\alpha_x f(x)=0
$$
in which all $\alpha_x\neq 0$. In particular, such an $f$ cannot vanish on all but one point of $X$.
\end{thm}

C. Prohata gave an elegant  proof of  Alon-F\"uredi's theorem as an application of Theorem \ref{Kar}. 

Let $AG(n, q)$ denote the $n$-dimensional affine space
over $\Fq$. 
A. Bishnoi, P. L. Clark, A. Potukuchi  and  J. R.  Schmitt gave the following nice application of Alon-F\"uredi's theorem in finite geometries in \cite{BCPS} Theorem 6.6..
\begin{thm}\label{BCPS_thm}
Let $S$ be a set of $k$ points in $AG(n, q)$. Then there are at least  $m(q,\ldots ,q, nq-k+1)-1$  hyperplanes of $AG(n, q)$ which do not meet $S$.
\end{thm}
Here $m(a_1, \ldots ,a_n; k)$ denotes the minimum value of the product $y_1\cdots y_n$, where $y_i$'s are integers  satisfying  $\sum_{i=1}^n y_i=k$ and $1\leq a_i\leq k$ for each $i$. 

A. Bishnoi used a special Gr\"obner basis technique, the footprint bound to prove   Alon-F\"uredi's theorem in his blogspot \cite{B}. This relates covering problems and the extensive Gr\"obner basis applications
for coding theory (see also \cite{B3}, \cite{C}). 

P. Sziklai  and Zs. Weiner proved the following nice generalization of  Alon-F\"uredi's theorem in \cite{SzW}. Their proof is based on M\"obius inversion and Zeilberger's method for proving binomial equalities.
\begin{thm}\label{Szw}
Let $0\leq k<n$ be integers.  Let $f\in I(\cV(n,k))$ and  suppose that $f(\ve v)\neq 0$ for each $\ve v\in \{0,1\}^n \setminus \cV(n,k)$. Then $\deg(f)>k$.  
\end{thm}

\begin{remark}
It is easy to verify that Theorem \ref{Szw} is sharp: consider $f(\ve x):=\prod_{j=0}^k (({\sum_{i=1}^n x_i)-j})\in \Q[\ve x]$. Then $f\in I(\cV(n,k))$ and  $f(\ve v)\neq 0$ for each $\ve v\in \{0,1\}^n \setminus \cV(n,k)$.  
\end{remark}

One of our main results is the following natural  generalization of  Sziklai--Weiner's Theorem. 

\begin{thm}\label{main}
Let $0\leq k<n$ be integers.  Let $f\in I(\cV(n,k))$ and  suppose that there exists a vector $\ve v\in \{0,1\}^n \setminus \cV(n,k)$ such that $f(\ve v)\neq 0$. Then $\deg(f)>k$.  
\end{thm}

We give here a very short proof for Theorem \ref{main}. Our simple proof is based on Gr\"obner basis theory.

Let $C_n$ denote the $n$-dimensional unit cube and $\Pi_{n-1}$  the $n$-dimensional   permutohedron. 

For a finite set of points $\cV \subseteq \F^n$ in an affine space and $\ve v\in \cV$ define
$\AC(\cV,\ve v)$ as the minimum size of an almost cover 
$\cV$ and $\ve v$. Define
$$
\AC(\cV):=\max_{\ve v\in \cV}\AC(\cV,\ve v)
$$  
and similarly
$$
\ac(\cV):=\min_{\ve v\in \cV}\AC(\cV,\ve v).
$$  

If $\cV$ is the vertex set of a convex polytope $P$, then we simply write $\ac(P)$. Thus, Alon--F\"uredi theorem  implies that $\ac(C_n)=n$ and we proved in \cite{HK} Theorem 1 that  $\ac(\Pi_{n-1})={n\choose 2}$.

We state our main results. To the best of the author's knowledge, our results are completely new. We combine Gr\"obner basis techniques with some simple upper bounds for binomial coefficients in our proofs.

First we give here a lower bound for the number $\AC(\cV)$ in terms of the size of $\cV$ and the dimension $n$. 

\begin{thm}\label{main2}
Let $\cV \subseteq {\F}^n$ be a finite subset of the affine space $\F^n$. 
Let $k\geq 0$ be any non-negative integer such that
$$
{n+k \choose n}< |\cV|.
$$
Then
$$
k<\AC(\cV),
$$
i.e. there exists a vector $\ve v\in \cV$ such that $k< AC(\cV,\ve v)$.
\end{thm}
\begin{remark}
It is easy to verify that Theorem \ref{main2} is sharp: let $q>1$ be a positive integer, $k:=q-1$ and $\cV:=J(n,q)$. Then $|\cV|={n+q-1 \choose  q-1}={n+k \choose n}$ and it follows from \cite{HR} Theorem 1.2 that $\AC(\cV)=k$.
\end{remark}

We can easily derive from Theorem \ref{main2} the following two Corollaries.

\begin{cor} \label{maincor1}
Let $\cV \subseteq {\F}^n$ be a finite subset of the affine space $\F^n$. Suppose that $|\cV|\geq 4^n$. Then 
$$
n<AC(\cV).
$$
\end{cor}

\begin{cor} \label{maincor2}
Let $\cV \subseteq {\F}^n$ be a non-empty finite subset of the affine space $\F^n$.  Then 
$$
\frac{n\cdot |\cV|^{\frac 1n}}{e}-n<AC(\cV),
$$
i.e. there exists a vector $\ve v\in \cV$ such that $\frac{n\cdot |\cV|^{\frac 1n}}{e} -n<AC(\cV,v)$.
\end{cor}

We prove a stronger version of Theorem \ref{main2} in the case of  $0-1$ vectors.

\begin{thm}\label{main3}
Let $\cV \subseteq \{0,1\}^n \subseteq {\F}^n$ be a finite $0-1$ subset.  
Let $0\leq k\leq n$ be non-negative integers such that
$$
\sum_{i=0}^k {n\choose i}<|\cV|.
$$
Then
$$
k<AC(\cV),
$$
i.e. there exists a vector $\ve v\in \cV$ such that $k< AC(\cV,\ve v)$.
\end{thm}
\begin{remark}
Clearly Theorem \ref{main3} is sharp: 
consider the  set $\cV:=\cV(n,k)$. Then it is easy to check that $|\cV|=\sum_{i=0}^k {n\choose i}$ and $AC(\cV)=k$. Namely the set of hyperplanes $\{(\sum_j x_j)-i:~ 1\leq i\leq k \}$ is an almost cover of the set $\cV(n,k)$.
\end{remark}

If the finite set of points $\cV$ has a transitive group of symmetries consisting of affine isomorphisms, then we can prove a lower bound for the number $\ac(\cV)$ in terms of the size of $\cV$ and the dimension $n$.

\begin{thm}\label{main4}
Let $\cV \subseteq {\F}^n$ be a finite subset of the affine space $\F^n$. Suppose  that there exists a finite subgroup $G\leq \Aff( {\F},n)$ such that $G$ is a transitive subgroup of the symmetric group $\Sym(\cV)$. Then the function $\AC(\cV,\cdot):\cV\to \N$ is constant, i.e., $\AC(\cV,\ve v)=\AC(\cV,\ve w)$ for each $\ve v\neq \ve w\in \cV$. 
\end{thm}

Finally Theorem \ref{main4} implies the following two Corollaries.

\begin{cor}\label{main5}
Let $\cV \subseteq {\F}^n$ be a finite subset of the affine space $\F^n$. Suppose  that there exists a finite subgroup $G\leq \Aff( {\F},n)$ such that $G$ is a transitive subgroup of the symmetric group $\Sym(\cV)$. Define $v:=|\cV|$. 
Let $k\geq 0$ be any non-negative integer such that
$$
{n+k \choose n}< v.
$$
Then
$$
k<\ac(\cV).
$$
\end{cor}

\begin{cor}\label{main6}
Let $\cV \subseteq {\F}^n$ be a non-empty finite subset of the affine space $\F^n$. Suppose  that there exists a finite subgroup $G\leq \Aff( {\F},n)$ such that $G$ is a transitive subgroup of the symmetric group $\Sym(\cV)$.   Then 
$$
\frac{n\cdot |\cV|^{\frac 1n}}{e}-n<ac(\cV).
$$
\end{cor}

\section{Preliminaries}






\subsection{Gr\"obner basis theory}

We recall some basic facts about  Gr\"obner 
bases   and standard monomials.  We refer to \cite{AL}, \cite{CLS} for details.
We say that a linear order $\prec$ on the monomials over 
variables $x_1,x_2,\ldots, x_n$ is a {\em term order}, if the following two conditions are satisfied:
\begin{itemize}
\item[(i)] 1 is
the minimal element of $\prec$;
\item[(ii)] $\ve u \ve w\prec \ve v\ve w$ holds for any monomials
$\ve u,\ve v,\ve w$ with $\ve u\prec \ve v$. 
\end{itemize}

The {\em leading monomial} ${\rm lm}(f)$
of a nonzero polynomial $f\in \F[x_1,x_2,\ldots, x_n]$ is the $\prec$-largest
monomial which appears with nonzero coefficient in the canonical form of $f$ as a linear
combination of monomials. Similarly, ${\rm lc}(f)$ denotes the leading coefficient of $f$, where $f\in \F[x_1,x_2,\ldots, x_n]$ is a nonzero polynomial.

Let $I$ be an ideal of the ring $\F[x_1,x_2,\ldots, x_n]$. Recall  that a finite subset $\cG\subseteq I$ is a {\it
Gr\"obner basis} of $I$ if for every $f\in I$ there exists a polynomial $g\in \cG$ such
that ${\rm lm}(g)$ divides ${\rm lm}(f)$. This means that the leading
monomials ${\rm lm}(g)$ for $g\in \cG $ generate the semi-group ideal of 
monomials
$\{ {\rm lm}(f):~f\in I\}$. It follows easily that  $\cG$ is
actually a basis of $I$, i.e. $\cG$ generates $I$ as an ideal of $\F[x_1,x_2,\ldots, x_n]$ (cf. \cite{CLS} Corollary 2.5.6). 
 It is a well--known result (cf. \cite[Chapter 1, Corollary
3.12]{CCS} or \cite[Corollary 1.6.5, Theorem 1.9.1]{AL}) that every
nonzero ideal $I$ of $\F[x_1,x_2,\ldots, x_n]$ has a Gr\"obner basis.

We say that a monomial $\ve w\in \F[\ve x]$ is a {\it standard monomial for $I$} if
it is not a leading monomial for any $f\in I$. We denote by  ${\rm Sm}(I)$ 
the set of standard monomials of $I$.
For a
nonzero ideal $I$ of $\F[\ve x]$ the set of monomials
${\rm Sm}(I)$ is a down-set: if $\ve w\in {\rm Sm}(I)$, $\ve u,\ve v$ are
monomials from $\F[\ve x]$ such that $\ve w=\ve u \ve v$ then 
$\ve u\in {\rm Sm}(I)$. 

If $\cV \subseteq {\F}^n$ is a finite subset of the affine space $\F^n$, then ${\rm Sm}(\cV)$ denotes the set of standard monomials of the ideal $I:=I(\cV)$.

We denote by $\func(\cV,\F)$ the $\F$-vector-space of all functions from $\cV$ to $\F$. It is a well-known that  ${\rm Sm}(\cV)$ constitutes a linear basis of the function space $\func(\cV,\F)$. Consequently $|{\rm Sm}(\cV)|=|\cV|$.

Let $\cF\subseteq 2^{[n]}$ be a set family. 
Then the characteristic vectors 
in $V(\cF)$ are all 0,1-vectors, hence the polynomials $x_i^2-x_i$ 
all vanish  on $V(\cF)$. We infer that the standard monomials of 
$I(\cF):=I(V(\cF))$ are square-free monomials.

Finally we introduce here shortly  the notion of reduction. Let $\prec$ be a fixed term order.
Let $\cG$ be a set of polynomials in $\F[x_1,\ldots,x_n]$ and
let $f\in \F[x_1,\ldots,x_n]$ be a fixed polynomial.
We can reduce $f$ by the set $\cG$ with respect to $\prec$.
This gives us a new polynomial $h \in \F[x_1,\ldots,x_n]$.

The term {\em reduction} means that we possibly repeatedly replace monomials
in $f$ by smaller ones (with respect to $\prec$).  This reduction process is the  following: if $w$ is a
monomial occurring in $f$ and ${\rm lm}(g)$ divides $w$ for some
$g\in {\cal G}$ (i.e. $w={\rm lm}(g)u$ for some monomial $u$), then
we replace $w$ in $f$ with $u({\rm lm}(g)-\frac{g}{lc(g)})$. It is easy to check that the
monomials in $u({\rm lm}(g)-\frac{g}{lc(g)})$
are $\prec$-smaller than $w$.

It is a basic fact that
${\rm Sm}(I)$ constitutes a basis of the $\F$-vector-space $\F[x_1,\ldots,x_n]/I$ in
the sense that  
every polynomial $g\in \F[x_1,\ldots,x_n]$ can be 
uniquely expressed as $h+f$ where $f\in I$ and
$h$ is a unique $\F$-linear combination of monomials from ${\rm Sm}(I)$. Consequently if $g\in \F[x_1,\ldots,x_n]$ is an arbitrary polynomial and $\cG$ is a Gr\"obner basis of $I$, then we can reduce $g$ with $\cG$ into a linear combination of
standard monomials for $I$.

In the case of the deglex term order a particular property of this reduced form $h$ is that $\mbox{deg}(h)\leq \mbox{deg}(g)$. Since our proofs rely on the well-known fact, that reduction modulo a a Gr\"obner basis with respect to a degree-compatible term order does not increase the total degree, we use the deglex term order in our proofs.

\subsection{Binomial coefficients}

We state here first 
the following well-known upper bound for the binomial coefficient (see \cite{BF} Exercise 4.2.1).

\begin{prop} \label{Bin_upper}
Let $1\leq k\leq n$ be integers.
Then
$$
{n \choose k} < \Big( \frac{ne}{k}\Big)^k.
$$
\end{prop}
\qed

The following Corollary is immediate. 

\begin{cor} \label{Bin_upper2}
Let $n\geq 1$, $k> -n$ be integers. Then
$$
{n+k \choose n}< e^n (1+ \frac kn)^n. 
$$
\end{cor}
\qed


\section{Proofs}

We use the following combinatorial Lemma in the proofs of our main results.

\begin{lemma} \label{monom}
Consider the set of monomials 
$$
D(n,k):=\{x^{\alpha}\in \F[x_1,\ldots,x_n]:~ \deg(x^{\alpha})\leq k\}.
$$
Then
$$
|D(n,k)|={n+k \choose n}.
$$
\end{lemma}

{\bf Proof of Theorem \ref{main}:}\\

It is easy to verify that 
$$
\Sm(\cV(n,k))=\{x_K:~ |K|\leq k\}.
$$
Namely each $x_L$, where $|L|>k$, is a leading monomial for the ideal $I(V(n,k))$.

Let $\ve v_T$ denote the vector with $f(\ve v_T)\ne 0$. Consider the  set $\cV:=\cV(n,k,T)$. Then 
$$
\Sm(\cV(n,k))\subseteq \Sm(\cV),
$$
hence
$$
\Sm(\cV(n,k))=\{x_K:~ |K|\leq k\}\subseteq \Sm(\cV).
$$
But $\Sm(\cV)$ is a down-set and  $|\Sm(\cV)|=|\cV|=|\cV(n,k)|+1$, which implies that there exists an $M\in { [n]\choose k+1}$ such that 
$$
\Sm(\cV)=\{x_K:~ |K|\leq k\}\cup \{x_M\}.
$$

We prove that $x_M(\ve v_T)\neq 0$. Namely, contrary suppose that $x_M(\ve v_T)=0$. Since $M\subseteq [n]$ and $|M|=k+1$, hence $x_M(\ve v_K)=0$ for each $K\in { [n]\choose {\leq k}}$.
Consequently $x_M\equiv 0$ in the space $\func(\cV,\F)$, a contradiction.

Clearly $x_M(\ve v_T)=1$, hence $x_M\equiv \chi_{\ve v_T}$ in the space $\func(\cV,\F)$.

Now let $f\in I(\cV(n,k))$ be a polynomial such that $c:=f(\ve v_T)\ne 0$. Then $f\equiv c\cdot\chi_{\ve v_T}\equiv c\cdot x_M$ in the space $\func(\cV,\F)$. 

Let $\cG$ denote a fixed deglex Gr\"obner basis of $I(\cV)$. Let $g$ denote the reduction of $f$ with $\cG$. Then it is easy to verify that  $g\equiv f \equiv c\chi_{\ve v_T}\equiv c x_M$ in the space   $\func(\cV,\F)$, and $x_M\in \Sm(\cV)$, consequently
$$
k+1=\deg(x_M)=\deg(g)\leq \deg(f). 
$$
\qed

{\bf Proof of Theorem \ref{main2}:}\\

Let $\cV \subseteq {\F}^n$ be a finite subset of the affine space $\F^n$. 
Let $k\geq 0$ be a non-negative integer such that
$$
|\cV|>{n+k \choose n}.
$$

First we prove that there exists a monomial $x^{\ve t}\in \Sm(\cV)$ such that $\deg(x^{\ve t})>k$.

Namely, if we  suppose contrary that 
$$
\Sm(\cV)\subseteq D(n,k),
$$
then
$$
|\cV|=|\Sm(\cV)|\leq |D(n,k)|={n+k \choose n},
$$
a contradiction.

Let $\ve w\in \cV$ be a fixed vector. Let  $\chi_{\ve w}\in \func(\cV,\F)$ denote the characteristic function of $\ve w$.  We can write up $\chi_{\ve w}$ as a sum of standard monomials of $I(\cV)$, since $\Sm(\cV)$ is a linear basis of the vector space $\func(\cV,\F)$:
$$
\chi_{\ve w}=\sum_{x^{\ve v}\in \Sm(\cV)} c_{\ve w,\ve v} x^{\ve v}.
$$

Since the set of characteristic vectors $\{\chi_{\ve z}: \ve z\in \cV\}$ is also a  basis of the space $\func(\cV,\F)$, there exists an $\ve u\in \cV$ vector such that the coefficient $c_{\ve u,\ve t}\ne 0$ in the expansion
$$
\chi_{\ve u}=\sum_{x^{\ve v}\in \Sm(\cV)} c_{\ve u,\ve v} x^{\ve v}.
$$

Let $H_1, \ldots ,H_m$ be an almost cover of $\cV$ and $\ve u$, which means that
$$
\cV\setminus \{\ve u \}\subseteq \cup_{j=1}^m H_j
$$
and $\ve u \notin \cup_{j=1}^m H_j$. Then there exist $L_j(\ve x)\in \F[\ve x]$ linear polynomials such that $H_j=V(L_j)$ for each $j$. Consider the polynomial 
$
F(\ve x):= (\prod_{j=1}^m L_j)(\ve x)
$. Let $\cG$ denote a fixed deglex Gr\"obner basis of the ideal  $I(\cV)$.  Let $g$ denote the reduction  of $F$ via $\cG$.  Then
$$
\deg(g)\leq \deg(F)=\deg(\prod_{j=1}^m L_j)= m.
$$
But there exists a $c\in \F$, $c\ne 0$ such that $g\equiv c\cdot \chi_{\ve u}$ in the space $\func(\cV,\F)$, because $H_1, \ldots ,H_m$ is an almost cover of $\cV$ and $\ve v$.  Then 
$$
g=c \cdot \sum_{x^{\ve v}\in \Sm(\cV)} c_{\ve u,\ve v} x^{\ve v},
$$
because $g$ is the reduction  of $F$ via $\cG$.

It follows that $\deg(g)>k$, since  $c_{\ve u,\ve t}\ne 0$ and $\deg(x^t)>k$. \qed 

{\bf Proof of Corollary \ref{maincor1}:}\\

This follows from the inequality
$$
{2n \choose n}<4^n.
$$

{\bf Proof of Corollary \ref{maincor2}:}\\

Let
$$
k:=\lfloor \frac{n\cdot |\cV|^{\frac 1n}}{e}-n \rfloor \in \Z.
$$
The assumption $|\cV|>0$ implies  that $k> -n$.  It follows from the definition of $k$ that 
$$
e^n (1+ \frac kn)^n\leq |\cV|. 
$$
If we combine this inequality with
$$
{n+k \choose n}< e^n \Big(1+ \frac kn\Big)^n
$$
appearing in Corollary \ref{Bin_upper2}, then we get 
$$
{n+k \choose n}<|\cV|.
$$
 Hence Corollary \ref{maincor2} follows from Theorem \ref{main2}. \qed

{\bf Proof of Theorem \ref{main3}:}\\

Let $\cV \subseteq \{0,1\}^n \subseteq {\F}^n$ be a finite set of $0-1$ vectors such that 
$$
|\cV|>\sum_{i=0}^k {n\choose i}.
$$
Then it is easy to verify that there exists a square-free monomial $x_T\in \Sm(\cV)$ such that 
$$
|T|=\deg(x_T)>k.
$$
Namely $\Sm(\cV)$ contains only square-free monomials. We argue in an indirect way. If $\Sm(\cV)\subseteq \{x_F:~ |F|\leq k\}$, then 
$$
|\cV|=|\Sm(\cV)|\leq \sum_{i=0}^k {n\choose i},
$$ 
a contradiction.

Let $\ve w\in \cV$ be a fixed vector. We can expand the  $\chi_{\ve w}\in \func(\cV,\F)$ characteristic function in the basis of standard   monomials as
$$
\chi_{\ve w}=\sum_{x_K\in \Sm(\cV)} c_{\ve w,K} x_K.
$$

Since the set of functions  $\{\chi_{\ve v}: \ve v\in \cV\}$ is a basis of the space $\func(\cV,\F)$,  there exists an $\ve y\in \cV$ vector such that the coefficient $c_{\ve y,T}\ne 0$ in the expansion 
$$
\chi_{\ve y}=\sum_{x_K\in \Sm(\cV)} c_{\ve y,K} x_K.
$$ 

Let $H_1, \ldots ,H_m$ be an almost cover of $\cV$ and $\ve y$ with affine hyperplanes, i.e.  
$$
\cV\setminus \{\ve y \}\subseteq \cup_{j=1}^m H_j
$$
and $\ve y \notin \cup_{j=1}^m H_j$. Let $k_i(\ve x)\in \F[\ve x]$ denote the linear polynomials such that $H_j=V(k_j)$ for each $j$. Define the polynomial 
$
h(\ve x):= (\prod_{j=1}^m k_j)(\ve x)
$. Let $\cH$ denote a fixed deglex Gr\"obner basis of the ideal  $I(\cV)$.
Let $d$ denote the reduction  of $h$ via the Gr\"obner basis $\cH$.  Then clearly 
$$
\deg(d)\leq \deg(h)\leq m.
$$

On the other hand  $h\equiv d\equiv d(\ve y)\cdot \chi_{\ve y}$ in the function space $\func(\cV,\F)$, where  $d(\ve y)\ne 0$.

Hence if we expand $d$ in the basis of standard monomials $\Sm(\cV)$, then the coefficient of $x_T$ in this expansion is $d(\ve y)\cdot c_{\ve y,T}\neq 0$. Consequently $\deg(d)>k$, because $|T|=\deg(x_T)>k$. \qed

{\bf Proof of Theorem \ref{main4}:}\\
Let $\ve v\neq \ve w\in \cV$ be arbitrary fixed vectors.

Let $H_1, \ldots ,H_m$ be a minimal set of affine hyperplanes, which is an almost cover of $\cV$ and $\ve v$, i.e. 
$$
\cV\setminus \{\ve v \}\subseteq \cup_{j=1}^m H_j
$$
and $\ve v \notin \cup_{j=1}^m H_j$. Then there exists an affine isomorphism  $g\in \Aff(\F,n)\cap \Sym(\cV)$ such that 
$g(\ve v)=\ve w$. This implies that  the affine hyperplanes $g(H_1), \ldots , g(H_m)$ consist of an almost cover of $\cV$, i.e. 
$$
\cV\setminus \{\ve w \}\subseteq \cup_j g(H_j)
$$
and  $g(\ve v)=\ve w \notin \cup_j g(H_j)$. This proves that 
$\ac(\cV,\ve v)=\ac(\cV,\ve w)$. \qed




\end{document}